\documentclass{amsart}
\usepackage{graphicx}
\usepackage{amssymb}
\usepackage{amsfonts}
\swapnumbers
\newtheorem{theorem}{Theorem}[section]

\newtheorem{proposition}[theorem]{Proposition}

\theoremstyle{definition}

\newtheorem{assumption}[theorem]{Assumption}

\numberwithin{equation}{section}
 \theoremstyle{plain}
 
 \numberwithin{equation}{section} 
 \numberwithin{figure}{section} 
 \theoremstyle{plain}
 \theoremstyle{remark}
 \newtheorem*{acknowledgement*}{Acknowledgement}

\newcommand{\ve}{{\varepsilon}}

\newcommand{\gam}{{\gamma}}

\newcommand{\bbR}{{\mathbb R}}
\newcommand{\bbS}{{\mathbb S}}

\begin{document}
\title[]{Large Deviations for Iterated Sums and Integrals}%
 \vskip 0.1cm
 \author{Yuri Kifer\quad\quad\quad\quad\quad\quad and \quad\quad\quad\quad\quad Ofer Zeitouni\\
Institute of Mathematics \quad\quad\quad\quad\quad\quad  Department  of Mathematics\\
 Hebrew University \quad\quad\quad\quad\quad\quad Weizmann Institute of Science\\
Jerusalem, Israel \quad\quad\quad\quad\quad\quad\quad\quad\quad\quad Rehovot, Israel}
\address{
Institute of Mathematics, The Hebrew University, Jerusalem 91904, Israel}
\email{ yuri.kifer@mail.huji.ac.il}%
\address{
Department of Mathematics, Weizmann Institute of Science, POB 26, Rehovot 761100, Israel}
\email{ofer.zeitouni@weizmann.ac.il}
\thanks{ }
\dedicatory{  }
 \date{\today}
\vskip 0.1cm

\begin{abstract}
We describe large deviations for normalized multiple iterated sums and integrals
of the form $\bbS_N^{(\nu)}(t)=N^{-\nu}\sum_{0\leq k_1<...<k_\nu\leq Nt}\xi(k_1)\otimes\cdots\otimes\xi(k_\nu)$,
$t\in[0,T]$ and $\bbS_N^{(\nu)}(t)=N^{-\nu}\int_{0\leq s_1\leq...\leq s_\nu\leq Nt}\xi(s_1)\otimes\cdots\otimes\xi(s_\nu)ds_1\cdots ds_\nu$, where $\{\xi(k)\}_{-\infty<k<\infty}$ and $\{\xi(s)\}_{-\infty<s<\infty}$ are centered bounded stationary vector processes whose sums or integrals satisfy
a trajectorial large deviations principle.
\end{abstract}
\maketitle

\markboth{ Yu.Kifer and O. Zeitouni}{Large Deviations}
\renewcommand{\theequation}{\arabic{section}.\arabic{equation}}
\pagenumbering{arabic}

\section{Introduction}\label{sec1}\setcounter{equation}{0}
Let $\{\xi(k)\}_{0\leq k<\infty}$ and $\{\xi(t)\}_{0\leq t<\infty}$ be discrete and continuous time $d$-dimensional
bounded stochastic processes. Namely, we assume throughout that $\sup_k|\xi(k)|,\,\sup_t|\xi(t)|\leq 1$ almost surely
though a uniform bound by any constant will do, as well..
We consider sequences of multiple iterated sums
\begin{equation}\label{1}
\bbS_n^{(\nu)}(t)=n^{-\nu}\sum_{0\leq k_1<...<k_\nu\leq tn}\xi(k_1)\otimes\cdots\otimes\xi(k_\nu),
\end{equation}
in the discrete time, and of multiple iterated integrals
\begin{equation}\label{2}
\bbS_n^{(\nu)}(t)=n^{-\nu}\int_{0\leq u_1<...<u_\nu\leq tn}\xi(u_1)\otimes\cdots\otimes\xi(u_\nu)du_1\cdots du_\nu
\end{equation}
where $t\in[0,T],\, T>0$. In the discrete time case we interpolate linearly to make $\bbS_n^{(\nu)}(t)$ Lipschitz
 continuous in $t$, while in the continuous time case $\bbS_n^{(\nu)}(t)$ is automatically Lipschitz continuous
 since $\xi$ is bounded, and so it is differentiable Lebesgue almost everywhere in $t$ (the Rademacher theorem)
 with a bounded derivative.
  Observe that we can write $\bbS_n^{(\nu)}(t)$  coordinate-wise in the following form
 \begin{equation}\label{3}
\bbS_n^{(\nu;i_1,...,i_\nu)}(t)=n^{-\nu}\sum_{0\leq k_1<...<k_\nu\leq tn}\xi_{i_1}(k_1)\cdots\xi_{i_\nu}(k_\nu)
\end{equation}
and
\begin{equation}\label{4}
\bbS_n^{(\nu;i_1,...,i_\nu)}(t)=n^{-\nu}\int_{0\leq u_1<...<u_\nu\leq tn}\xi_{i_1}(u_1)\cdots\xi_{i_\nu}(u_\nu)du_1
\cdots du_\nu
\end{equation}
 We will be interested in large deviations of $\bbS_n^{(\nu)}$ as $n\to\infty$.
 
  The sequences of multiple iterated sums and of multiple iterated integrals
 were called signatures in recent papers related to the rough paths theory, data sciences and machine
learning (see, for instance, \cite{HL}, \cite{DR}, \cite{DET} and references there). 
  Among important applications we have in mind is the case of dynamical systems when
 we have stationary processes $\xi(k)=g\circ f^k$ or $\xi(t)=g\circ f^t$ where $f^k$ or $f^t$ are  discrete or 
 continuous time measure preserving transformations on, say, a compact space and $g$ is a continuous function, so that our
 boundedness assumption on processes $\xi$ will hold authomatically.
 Observe that several limit theorems for iterated sums and integrals were obtained 
 in \cite{Ki23} and \cite{Ki24}, so the present large deviations result complements the picture.

\section{Main results}\label{sec2}\setcounter{equation}{0}
We will proceed first with the continuous time case while the discrete time case after the linear interpolation can
be treated essentially similar. Set
\[
\phi_n(t)=n^{-1}\int_0^{tn}\xi(u)du.
\]
Then
\[
\xi(s)=\frac {d\phi_n(s)}{ds}=\dot\phi_n(s).
\]
The distribution of the process $\xi$ induces a law of $\phi_n(t),\, t\in [0,T]$ which we call $\mu_n$ which does 
not lead to a confusion. Our results will be obtained assuming that the trajectorial large deviations principle (LDP)
holds true for $\phi_n$ which is a standard problem. To state it, let $H$ be the space of Lipschitz continuous paths 
$\gam(t),\, t\in[0,T]$ in $\bbR^d$ with $\gam(0)=0$, having the Lipschitz constant 1 (any other constant $C>0$ in place
of 1 will do, as well) and equiped with the supremum norm. Observe that $H$ is a compact subset of $C([0,T];\bbR^d)$
and any element of $H$ is Lebesgue almost everywhere differentiable on $[0,T]$ by the Rademacher theorem.

\begin{assumption}\label{as} The random vectors $\xi$ are bounded ($\sup_{0\leq t<\infty}|\xi(t)|\leq 1$) and 
the sequence of paths $\phi_n$ satisfies LDP in $H$ with a rate function $I$
(which is automatically good since $H$ is compact) meaning that for any open set $G\subset H$ and closed set
$F\subset H$,
\begin{eqnarray}\label{5}
&\liminf_{n\to\infty}\frac 1n\log\mu_n\{G\}\geq-\inf_{\gam\in G}I(\gam)\,\,\,\mbox{and}\\
&\limsup_{n\to\infty}\frac 1n\log\mu_n\{F\}\leq-\inf_{\gam\in F}I(\gam).\nonumber
\end{eqnarray}
\end{assumption}
We remark that in Assumption  \ref{as}, one could replace the condition $\sup_{0\leq t<\infty}|\xi(t)|\leq 1$ by
$\sup_{0\leq t<\infty}|\xi(t)|\leq C$ with any constant $C$, with minimal changes in the arguments.

Consider the map $\Phi^{(\nu)}:\, H\to C([0,T];\,\bbR^{d\nu})$ from $H$ to the space $C([0,T];\,\bbR^{d\nu})$ of continuous
 $d\nu$-dimensional vector functions on $[0,T]$, acting by
\[
\Phi^{(\nu)}(\gam)(s)=\int_{0\leq u_1\leq\cdots\leq u_\nu\leq s}\dot\gam(u_1)\otimes\cdots\otimes
\dot\gam(u_\nu)du_1\cdots du_\nu,\,\, s\in[0.T].
\]
In particular, we have
\[
\Phi^{(\nu)}(\phi_n)(t)=\int_{0\leq u_1<...<u_\nu\leq t}\xi(nu_1)\otimes\cdots\otimes\xi(nu_\nu)du_1\cdots du_\nu=\bbS_n^{(\nu)}(t).
\]

Our main result here is the following theorem.
\begin{theorem}\label{th}
Let the Assumption \ref{as} hold true. Then $S_n^{(\nu)}$ satisfies LDP in $C([0,T];\,\bbR^{d\nu})$ with the good rate function
\[
\tilde I(\psi)=\inf\{ I(\gam):\,\gam\in H,\, \psi=\Phi^{(\nu)}(\gam)\}.
\]
\end{theorem}

The main step in the proof of Theorem \ref{th} is the following proposition where we show that the map $\Phi^{(\nu)}$ is
 H\" older continuous with respect to the supremum norm. 

\begin{proposition}\label{pr} The map $\Phi^{(\nu)}(\gam)$ is $\frac 12$-H\"older continuous in 
$\gam$ as map from $H$ to  $C([0,T];\bbR^{d\nu})$, and further 
 $\Phi^{(\nu)}(\gam)(t)$ is Lipschitz continuous in $t$ whenever $\gamma\in H$.
If we assume that $\dot\gam(t),\, t\in[0,T]$ is bounded almost surely by some constant $C>0$ and not by 1, 
then the claim remains true with some other Lipschitz constant.
\end{proposition}
Before proving this claim we will show first its application to our large deviations problem.
By Assumption \ref{as} the LDP for the integrals $\eta_n(t)=n^{-1}\int_0^{tn}\dot\gam(u)du$ holds true
with a good rate function  and $\Phi^{(\nu)}$  is a continuous map, and so
 the contraction principle (see, for instance, Theorem 4.2.1 in \cite{DZ}) can be applied here. It yields that
\[
\tilde I(y)=\inf\{ I(\gam):\,\gam\in H,\, y=\Phi^{(\nu)}(\gam)\}
\]
is a good rate function for $\bbS_n^{(\nu)}$ and Theorem \ref{th} follows. Moreover, we obtain from here the LDP 
coordinate wise by the contraction principle.
 Some large deviations results for iterated sums as above in the one dimensional case $d=1$ were obtained in the
  $U$-statistics research field (see \cite{EL}) but only when the discrete time processes $\{\xi(k)\}_{0\leq k<\infty}$
   appearing above consist of i.i.d. random variables.

\section{Proofs}\label{sec3}\setcounter{equation}{0}
Next, we will prove Proposition \ref{pr}.
\begin{proof} 
We will proceed by induction in $\nu$. For $\nu=1$ the assertion clearly holds true and $\Phi^{(1)}_n$ is, in
fact, Lipschitz continuous in all arguments with the constant 1. Assume that Assertion holds true for $\Phi^{(\nu)}$
with $\nu=1,...,k$ and prove it for $\Phi^{(k+1)}$. We have
\begin{equation}\label{6}
\Phi^{(k+1)}(\gam)(t)=\int_0^{t}\dot\gam(u)\otimes\Phi^{(k)}(\gam)(u)du.
\end{equation}
Since we assume that $|\dot\gam(u)|\leq 1,\, u\in[0,T]$ a.s. we obtain easily that 
$|\Phi^{(k)}(\gam)(u)|\leq \frac {T^k}{k!}$ a.s. in $u\in[0,T]$, where 
$\Phi^{(k)}(\gam)(u)$ is considered as a vector in $\bbR^{dk}$ and $|\cdot|$
is a vector norm. It follows that $\Phi^{(k+1)}(\gam)(t)$ is Lipschitz continuous in the time parameter $t$.

Now, let $l\ve\leq u<(l+1)\ve$ where $\ve>0$. Then, taking into account that $\Phi^{(k)}(\gam)(u)$ is Lipschitz continuous
in $u$ by the induction hypothesis we have
\begin{equation}\label{7}
|\int_{l\ve}^{(l+1)\ve}(\dot\gam(u)\otimes\Phi^{(k)}(\gam)(u)-\dot\gam(u)\otimes\Phi^{(k)}(\gam)(l\ve))du |\leq C\ve^2
\end{equation}
for some $C>0$.

Set $u(\ve)=l\ve$ if $l\ve\leq u<(l+1)\ve$. Let $\gam$ and $\tilde\gam$ satisfy
\begin{equation}\label{8}
|\gam(t)-\tilde\gam(t)|\leq\ve^2\quad\mbox{a.s. in $t\in[0,T]$}.
\end{equation}
Observe that for any $t>0$,
\begin{equation}\label{9}
\int_0^t\dot\gam(u)\Phi^{(k)}(\gam(u(\ve))du=\int_0^t(\gam(\min(t,u(\ve)+\ve))-\gam(u(\ve)))\Phi^{(k)}(\gam)(u(\ve))du
\end{equation}
and the same holds true for $\tilde\gam$ in place of $\gam$. This together with the induction hypothesis yields that
\begin{equation}\label{10}
\sup_{0\leq t\leq T}|\int_0^{t}\big(\dot\gam(u)\Phi^{(k)}(\gam)(u(\ve))-\dot{\tilde\gam}(u)
\Phi^{(k)}(\tilde\gam)(u(\ve))\big)du|\leq\tilde C\ve
\end{equation}
for some constant $\tilde C>0$ where we used $[t/\ve]$ times (\ref{7}). Now Proposition \ref{pr} follows from (\ref{1}), (\ref{2}), (\ref{4})  and (\ref{5}).

\end{proof}

\section{Law of large numbers}\label{sec4}\setcounter{equation}{0}

We observe that the map $\Phi^{(\nu)}$ proved above to be continuous can be used also to obtain certain law of 
large numbers result for iterated sums and integrals. Namely, suppose that 
\begin{equation}\label{11}
\lim_{n\to\infty}\phi_n(t)=tQ \quad\quad\mbox{a.s.},
\end{equation}
i.e. the standard law of large numbers for sums $\sum_{0\leq k\leq th}\xi(k)$ or integrals $\int_0^{tn}\xi(u)du$ 
holds true. Since $\Phi^{(\nu)}$ is continuous we obtain that almost surely,
\begin{eqnarray}\label{12}
&\lim_{n\to\infty}\bbS_n^{(\nu)}(t)=\lim_{n\to\infty}\Phi^{(\nu)}(\phi_n)(t)=\Phi^{(\nu)}(Q)(t)\\
&=Q^{\otimes\nu}\int_{0\leq u_1\leq...\leq u_\nu\leq t}du_1\cdots du_\nu=\frac {Q^{\otimes\nu} t^\nu}{\nu!}\nonumber
\end{eqnarray}
where $Q^{\otimes\nu}=Q\otimes\cdots\otimes Q$ is the $\nu$ times tensor product of $Q$ which is the $d\nu$-dimensional
vector with coordinates $Q^{\otimes\nu}_{i_1,...,i_\nu}=Q_{i_1}\cdots Q_{i_\nu}$. 

We observe that the continuity of the map $\Phi^{(\nu)}$ was proved here under the assumption that $\xi(k),\,\xi(t)$
are uniformly bounded processes, and so the law of large numbers (\ref{2}) is obtained here under this assumption which is
a restriction. In \cite{Ki25} this law was established by another method under substantially more general conditions, 
essentially the same as in the standard $\nu=1$ case.



\begin{thebibliography}{99}

\bibitem{BBK} P. Baldi, G. Ben Arous, G. Kerkyacharian, {\em Large deviations and the Strassen theorem in
H\" older norm}, Stoch. Proc. Appl. 42 (1992), 171--180.

















\bibitem{DZ} A. Dembo and O. Zeitouni, {Large Deviations Techniques and Applications}, Springer 2nd. ed., New York, 1998.







\bibitem{DET} J. Diehl, K. Ebrahimi-Fard and Nicolas Tapia, {\em Generalized iterated-sums signatores},
arXiv:2012.04597



\bibitem{DR} J. Diehl and J. Reizenstein, {\em Invariants of multidimensional time series based
 on their iterated-integral signature}, Acta Appl. Math. 164 (2019), 83--122.




\bibitem{EL} P. Eichelshacher and M. L\" owe, {\em A Large Deviation Principle for $m$-Variate
 von Mises-Statistics and U-Statistics}, J. Theoretical Probab. 8 (1995), 807--824.
























 \bibitem{HL} B. Hambly and T. Lyons, {\em Uniqueness for the signature of a path of bounded variation
 and the reduced path group}, Ann. Math. 171 (2010), 109--167.









\bibitem{Ki23} Yu. Kifer, {\em Limit theorems for signatures}, arXiv: 2306.13376.


\bibitem{Ki24} Yu. Kifer, {\em Almost sure approximations and laws of iterated logarithm for
 signatures}, Stoch. Proc. Appl. 182 (2025), 104576.



\bibitem{Ki25} Yu. Kifer, {\em Iterated ergodic theorems}, arxiv: 2501.15633













\end{thebibliography}


\end{document}